\documentclass[11pt]{article}

\usepackage[utf8]{inputenc}

\usepackage{epsfig}
\usepackage{amsmath,amsfonts,amssymb}
\usepackage{color}
\definecolor{vertfonce}{rgb}{0.20, 0.46, 0.25}
\definecolor{rougefonce}{rgb}{0.64, 0.09, 0.20}
\usepackage[ps2pdf,
            breaklinks=true,
            colorlinks=true,
            linkcolor=rougefonce,
            citecolor=vertfonce]{hyperref}

\input{macro}

\title{Symplectic inverse spectral theory for pseudodifferential
  operators}%
\author{V\~{u} Ng\d oc San\footnote{IRMAR (UMR 6625), Universit{\'e} de
    Rennes 1, Campus de Beaulieu, 35042 Rennes cedex (France)}}
\date {June 2008}
\begin{document}

\maketitle
\begin{abstract}
  We prove, under some generic assumptions, that the semiclassical
  spectrum modulo $\mathcal{O}(\h^2)$ of a one dimensional
  pseudodifferential operator completely determines the symplectic
  geometry of the underlying classical system. In particular, the
  spectrum determines the hamiltonian dynamics of the principal
  symbol.
\end{abstract}

\section{Introduction}

In this article I would like to advocate an inverse spectral theory
for pseudodifferential operators. What does this means ?  One of the
most famous inverse spectral problems, made fashionable by Kac's very
entertaining article~\cite{kac}, with a mind-catching title ``Can one
hear the shape of a drum ?''\footnote{Kac attributes the problem to
  Bochner and the title to Bers.}, was about the Laplace operator on a
bounded domain $\Omega\subset \RM^n$. Frequencies $\nu$ solutions to
the eigenvalue problem
\[
\frac12\Delta u = \nu^2 u, \qquad u=0 \text{ on } \partial\Omega
\]
may be viewed as harmonics that can be heard when the interior of the
``membrane'' $\Omega$ vibrates freely. The question was whether the
knowledge of \emph{all} frequencies completely determines $\Omega$ (up
to isometry, of course). As Kac mentioned, this question appears
naturally in the context of quantum mechanics, for a particle trapped
in a hard potential well. An important observation in this paper was
the relevance of the \emph{Weyl law}, which let us find the
\emph{volume} (or area when $n=2$) of $\Omega$ from the
\emph{asymptotic behaviour of large eigenvalues}.

Counterexamples are now known~: there are non-isometric shapes in
$\RM^2$ that produce different
frequencies~\cite{gordon-webb-wolpert}. Nonetheless, this fact should
not let us think that the problem has become obsolete. As far as I
know, the seemingly simple case of a convex, bounded domain
$\Omega\in\RM^2$ with analytic boundary is still open (although by the
time this article is published, it might very well have been settled
by Zelditch, see~\cite{zelditch-inverse-I,zelditch-inverse-survey}).

Understanding this problem requires putting it in a wider
perspective. A natural variant of Kac's problem is whether the
spectrum of the Laplace operator $\Delta_g$ on a compact riemannian
manifold $(M,g)$ determines the metric $g$. Although, here again,
counterexamples have been known for a long time~\cite{milnor-laplace},
our understanding remains relatively poor. Recent works by Zelditch
and Guillemin suggest that microlocal tools are quite relevant for all
these questions. This, in turn, is a hint that more general operators
than the Laplacian could be dealt with similarly.

From a quantum mechanical viewpoint, Kac's situation is quite
extreme. A more natural setting would involve a particle `trapped' by
a smooth potential well. No more boundary problems, but instead a
Schrödinger operator on $\RM^n$
\[
P = -\frac{\h^2}2\Delta + V(x).
\]
Of course now the \emph{potential function} $V$ should be recovered
from the spectrum of $P$. This inverse spectral problem has been
studied a lot, but only very recently have microlocal tools similar to
those used by Guillemin and Zelditch been applied to
it~\cite{guillemin-uribe-inverse,colin-guillemin}.

Here, I would like to shift again the initial problem one step further
away. Instead of the Laplacian, or the Schrödinger operator, why not
consider \emph{any} (elliptic) differential operator, or even, since
we're at it, any \emph{pseudodifferential} operator ? Of course, since
there is no domain $\Omega$ anymore, no potential function $V$, the
sensible question is \emph{what} should we try to recover from the
spectrum ?

The inverse spectral problems I've mentioned here can all be
understood as \emph{semiclassical limits}. From a quantum object, the
spectrum, one wants to recover classical observables such as the
metric $g$, or the potential $V$. These quantities, in turn, fully
determine the \emph{classical dynamics} of the system. For general
pseudodifferential operators, semiclassical analysis still shows the
strong relationship between the classical dynamics and the quantum
spectrum, so I believe that the most natural ``object'' that we should
try and recover from the spectrum is precisely this classical
dynamics. This, precisely, amounts to determining the \emph{principal
  symbol} of the operator. In fact, if we keep in mind Weyl's
asymptotics, this sounds fairly natural, for it is well known that
Weyl's asymptotics extend to arbitrary pseudodifferential operators,
provided that we compute phase space volumes defined by energy ranges
given by the principal symbol~\cite{robert-weyl}.

As in the riemannian case, one should take into account a symmetry
group acting on the classical data. For general pseudodifferential
operators, there's only one available~: the group of
symplectomorphisms, acting on the phase space $M$. This is a much
bigger group than the group of riemannian isometries, in accordance
with the fact that the space of principal symbols $\Cinf(M)$ is much
bigger than the space of riemannian metrics, or potential functions.

\section{The setting}

Since we aim at recovering the classical dynamics from the spectrum,
we are going to work in the setting of semiclassical
pseudodifferential operators, which we recall here. Throughout this
work, we only consider the \emph{one-dimensional} theory. It would be
very interesting to have higher dimensional results, but it is not
expected that such precise results would persist. However, a
reasonable challenge would be to undertake a similar study for the
\emph{completely integrable} case.

The classes $\Psi^d(m)$ of semiclassical pseudodifferential operators
we use are standard.  Let $M=T^*\RM=\RM^2_{(x,\xi)}$. Let $d$ and $m$
be real numbers. Let $S^d(m)$ be the set of all families
$(p(\cdot;\h))_{\h\in (0,1]}$ of functions in $\Cinf(M)$ such that
\begin{equation}
  \forall \alpha\in\NM^2, \quad \abs{\partial^\alpha_{(x,\xi)}
    p(x,\xi;\h)} \leq C_\alpha\h^d(1+\abs{x}^2+\abs{\xi}^2)^{\frac{m}2},
  \label{equ:symboles}
\end{equation}
for some constant $C_\alpha>0$, uniformly in $\h$.  Then $\Psi^d(m)$
is the set of all (unbounded) linear operators $P$ on $L^2(\RM)$ that
are $\h$-Weyl quantisations of symbols $p\in S^d(m)$~:
\[
(Pu)(x) = (Op^w_\h(p)u)(x) =
\frac{1}{2\pi\h}\int_{\RM^{2}}\!\!\!e^{\frac{i}{\h}\pscal{x-y}\xi}
p({\textstyle\frac{x+y}2},\xi;\h) u(y)\abs{dyd\xi}.
\]
The number $d$ in~(\ref{equ:symboles}) is called the $\h$-order of the
operator. Unless specified, it will always be zero here.  In this work
all symbols are assumed to admit a ``classical'' asymptotic expansions
in integral powers of $\h$ (that is to say, in the ladder
$(S^d(m))_{d\in\ZM, d\geq d_0}$ for some $d_0\in\ZM$). The leading
term in this expansion is called the principal symbol of the operator.

Thus, the Schr{\"o}dinger operator $P=-\frac{\h^2}2\Delta + V$ on
$\RM$ is a good candidate, of $\h$-order zero, whenever $V$ has at
most a polynomial growth.

We use in this article the standard
properties of such pseudodifferential operators.  In particular the
composition sends $\Psi^d(m)\times\Psi^{d'}(m')$ to
$\Psi^{d+d'}(m+m')$. Moreover all $P\in\Psi^0(0)$ are bounded:
$L^2(\RM)\fleche L^2(\RM)$, uniformly for $0<\h\leq 1$.

An operator $P\in\Psi(m)$ is said to be \emph{elliptic at infinity} if
there exists a constant $C>0$ such that the principal symbol $p$
satisfies
\[
\abs{p(x,\xi)}\geq \frac1C(\abs{x}^2+\abs{\xi}^2)^{m/2}
\]
for $\abs{x}^2+\abs{\xi}^2\geq C$.

If $P$ has a real-valued Weyl symbol, then it is a symmetric operator
on $L^2$ with domain $\Cinf_0(\RM)$. If furthermore the principal
symbol is elliptic at infinity, then $P$ is essentially selfadjoint
(see for instance~\cite[proposition 8.5]{dimassi-sjostrand}).

Finally, when $P\in\Psi^0(m)$ is selfadjoint and elliptic at infinity,
then for any $f\in\Cinf_0(\RM)$, the operator $f(P)$ defined by
functional calculus satisfies
$f(P)\in\cap_{k\in\NM}(\Psi^0(-km))$. See for
instance~\cite{dimassi-sjostrand} or~\cite{robert} for details.

The advantage of the semiclassical theory is that it allows us to use
richer versions of Weyl's asymptotics. Instead of considering the
limit of large eigenvalues, we fix a bounded \emph{spectral window}
$I=[E_0,E_1]\subset \RM$ and study the asymptotics of all eigenvalues
in $I$, as $\h\fleche 0$.
\begin{defi}
  \label{defi:hypo}
  We say that Assumption $\mathcal{A}(P,\mathcal{J},I)$ holds
  whenever
  \begin{enumerate}
  \item $P$ is a selfadjoint pseudodifferential operator in
    $\Psi^0(m)$ with principal symbol $p$, elliptic at infinity.
  \item $\mathcal{J}\subset[0,1]$ is an infinite subset with zero as
    an accumulation point.
  \item There exists a neighbourhood $J$ of $I$ such that $p^{-1}(J)$
    is compact in $M$.
  \end{enumerate}
\end{defi}
If Assumption $\mathcal{A}(P,\mathcal{J},I)$ holds, we denote by
$\Sigma_\h(P,I)$ the spectrum of $P=P(\h)$ in $I$ (including
multiplicities). We denote by $\Sigma(P,\mathcal{J},I)$ the family of
all $\{\Sigma_\h(P,I); \quad \h\in\mathcal{J}\}$.  It is well known
that $\Sigma_\h(P,I)$ is discrete for $\h$ small enough (see
eg.~\cite[Théorème 3.13]{robert}, in a slightly different
setup). Notice that when $m>0$, the properness condition 3. is always
satisfied.
\begin{prop}
  Let $P$ be a selfadjoint pseudodifferential operator in $\Psi(m)$,
  with principal symbol $p$, elliptic at infinity. Let $J\subset \RM$
  be a closed interval such that $p^{-1}(J)$ is compact. Then for any
  open interval $I\subset J$ there exists $\h_0>0$ such that the
  spectrum of $P$ in $I$ is discrete for $\h\leq \h_0$.
\end{prop}
\begin{demo} The case $m>0$ is probably the most standard. We recall
  it quickly.
  \paragraph{Case $m>0$. ---} Let $\chi\in\Cinf_0(J)$ be equal to 1 on
  $I$. Then by pseudodifferential functional calculus, $f(P)$ is
  compact for $\h$ small enough. Therefore, denoting by $\Pi_I$ the
  spectral projector on $I$, we have that $\Pi_I=\Pi_I f(P)$ is
  compact. This implies that $\Pi_I$ has finite rank~: the spectrum in
  $I$ is discrete.

  \paragraph{Case $m\leq 0.$ ---}
  First we show that one can replace $J$ by an unbounded interval
  containing $I$. Thus assume $J$ is compact.  For notational
  convenience we let $J=[-1,0]$. For any $\alpha\in\,]{-1},0[$, Sard's
  theorem ensures the existence of a regular value $\lambda\in\,
  ]\alpha,0[$ for $p$. Then $\mathcal{C}:=p^{-1}(\lambda)$ is a
  compact 1-dimensional submanifold of $\RM^2$~: it is a finite union
  of circles. Let $\Omega$ be the unbounded component of
  $\RM^2\setminus \mathcal{C}$. Suppose first that $p_{\restr
    \Omega}>\lambda$. Since the differential of $p$ does not vanish on
  $\mathcal{C}$, $p<\lambda$ in all the bounded components of
  $\RM^2\setminus \mathcal{C}$. Therefore
  $p^{-1}(\,]{-\infty},\lambda])$ is compact and one may replace $J$
  by $]{-\infty},\lambda]$. Now if on the contrary $p_{\restr
    \Omega}<\lambda$ on the unbounded components, we have to apply the
  same argument for $\lambda'\in\, ]{-1},\alpha[$. Because
  $\lambda'<\lambda$, the new bounded components contain the old one,
  and therefore consist of the points where $p>\lambda'$. Then
  $p^{-1}([\lambda',+\infty[\,)$ is compact, and one may replace $J$
  by $[\lambda',+\infty[$.

  For the rest of the proof, we may suppose that $J=\,]{-\infty},
  0]$. Let $-\epsilon\in J\setminus I$, close to the origin. Let
  $\chi\in\Cinf(\RM)$ be such that
  \[
  \begin{cases}
    \chi(x)=  -\epsilon/2 & \text{ for } x\in\,]{-\infty},-\epsilon]\\
    \chi(x) = x & \text{  for } x\geq -\epsilon/3\\
    \chi(x) \geq -\epsilon/2 & \text{ everywhere.}
  \end{cases}
  \]
  Let $p_{\h}$ be the Weyl symbol of $P$ and define
  $\tilde{p}_{\h}=\chi\circ p_{\h}$. Then $\tilde{p}_{\h}$ is a symbol
  in $\Psi(m)$, with $\tilde{p}_{\h}\geq -\epsilon/2$ on $\RM^2$ and
  $\tilde{p}_{\h}=p_{\h}$ outside the set
  $p_{\h}^{-1}(\,]{-\infty},-\epsilon/3])$. Because $m\leq 0$, the set
  $p_{\h}^{-1}(\,]{-\infty},-\epsilon/3])$ is included in $p^{-1}(J)$
  for $\h$ small enough and hence must be compact as well.  Denote by
  $\tilde{P}$ the Weyl quantisation of $\tilde{p}$. Then for $\h$
  small enough $(\tilde{P}-\lambda)$ is invertible for all $\lambda\in
  I$. We can write
  \[
  P-\lambda =
  (\tilde{P}-\lambda)\left(\textup{Id}+(\tilde{P}-\lambda)^{-1}(P-\tilde{P})\right).
  \]
  Since $P-\tilde{P}$ is the Weyl quantisation of a compactly
  supported symbol (with support in a compact independent of $\h$), it
  is of the trace class. By analytic Fredholm theory, we may take the
  determinant of $(\textup{Id}+(\tilde{P}-\lambda)^{-1}(P-\tilde{P}))$
  and conclude that the spectrum of $P$ consists of the zeroes (with
  multiplicities) of a non-vanishing holomorphic function and hence is
  discrete.
\end{demo}

The goal of this article is to recover the dynamics of the hamiltonian
$p$ in the region $p^{-1}(I)$ for any operator $P$ for which
Assumption $\mathcal{A}(P,\mathcal{J},I)$ holds, for some subset
$\mathcal{J}\subset[0,1]$.  Of course if we can do it for an arbitrary
compact interval $I\subset\RM$, we recover the full dynamics of $p$.

It turns out that, under some genericity conditions, these inverse
spectral problems are fairly easy, compared to the general
multi-dimensional problems alluded to in the introduction, in the
sense that they only require a few terms in the asymptotics of the
spectrum. Having this in mind, for $\alpha\in\RM$ we denote by
$\Sigma(P,\mathcal{J},I) + \mathcal{O}(\h^\alpha)$ the equivalence
class of all $\Sigma_\h(P,I)$ modulo $\h^\alpha$. Our main result is
Theorem~\ref{theo:sympl}, but we also state several intermediate results
that require weaker hypothesis. An informal statement of
Theorem~\ref{theo:sympl} is as follows.
\begin{theo}[Theorem~\ref{theo:sympl}]
  Let Assumption $\mathcal{A}(P,\mathcal{J},I)$ hold, and denote
  $M=p^{-1}(I)$. Suppose that $p_{\restr M}$ is a Morse function.
  Assume that the graphs of the periods of all trajectories of the
  hamiltonian flow defined by $p_{\restr M} $, as functions of the
  energy, intersect generically.

  Then the knowledge of $\Sigma(P,\mathcal{J},I)+\mathcal{O}(\h^2)$
  determines the dynamics of the hamiltonian system $p_{\restr M}$.
\end{theo}
In fact, we determine completely the Hamiltonian $p$ up to symplectic
equivalence. Perhaps the most difficult step, for which Weyl's
asymptotics are not enough, is the seemingly simple question to count
the number of connected components of $p^{-1}(E)$, for a regular
energy $E\in I$ (Theorem~\ref{theo:number}).

\ouf

Although we state everything for pseudodifferential operators defined
on on $\RM$, it is most probable that all results extend to the case
of pseudodifferential operators defined on a one-dimensional compact
manifold equipped with a smooth density, and to the case of Toeplitz
operators on two-dimensional symplectic manifolds.

\ouf

The plan of the paper follows a fairly logical progression.  Since we
always work modulo symplectomorphisms, it is not reasonable to look
for a formula that would give the principal symbol $p$. Instead we
will try to recover as many \emph{symplectic invariants} as possible
from the spectrum so that, given two spectra, we should be able to
tell whether they come from isomorphic systems.

Thus, the geometric object under study is a proper map
$p:M\fleche\RM$, where $M$ is a symplectic 2-manifold.  The simplest
symplectic invariants of this map are in fact topological invariants,
and are dealt with in Sections~\ref{sec:sing}
and~\ref{sec:topo}. Indeed, it follows from the action-angle theorem
that as soon as $E\in\RM$ is a regular value of $p$, then the fibres
of $p$ consist of a finite number of closed loops, each one
diffeomorphic to a circle. Therefore, we need to be able to detect
\begin{enumerate}
\item Whether an energy $E\in\RM$ is a regular or critical value of
  $p$; this is done in Section~\ref{sec:sing}
  (Theorem~\ref{theo:sing}).
\item When $E$ is a regular value, the number of connected components
  of the fibre $p^{-1}(E)$; Section~\ref{sec:connected} discusses this
  point (Theorem~\ref{theo:number}).
\end{enumerate}
Putting these results together we are able to recover the topological
type of the singular fibration (Theorem~\ref{theo:topo}). Then in
Section~\ref{sec:sympl}, relying on the classification result of
Dufour-Molino-Toulet~\cite{dufour-mol-toul,toulet-these} (and some
additional argument) we finally manage to recover the symplectic
geometry of the system (Theorem~\ref{theo:sympl}). 

\section{Singularities}
\label{sec:sing}

In order to detect whether a given energy $E_0\in\RM$ is a critical
value of $p$ or not, it is enough to know the spectrum of $P$ in a
small ball around $E_0$, at least under some nondegeneracy conditions.

Recall that a function $f:M\fleche\RM$ is said to have a nondegenerate
critical point $m\in M$ when $df(m)=0$ and the Hessian $f''(m)$ is a
nondegenerate quadratic form. Since $M$ has dimension 2, there are
only two cases~:
\begin{enumerate}
\item \textbf{Elliptic case~:} there are local symplectic coordinates
  $(x,\xi)$ in $T_mM$ such that $f''(m)(x,\xi)=C(x^2+\xi^2)$, for some
  constant $C\neq 0$.
\item \textbf{Hyperbolic case~:} there are local symplectic
  coordinates $(x,\xi)$ in $T_mM$ such that $f''(m)(x,\xi)=Cx\xi$, for
  some constant $C\neq 0$.
\end{enumerate}
We refer to each of these two cases as the \emph{type} of the
singularity $m$.
\begin{theo}
  \label{theo:sing}
  Let $I$ be an interval containing $E_0$ in its interior, and let
  Assumption $\mathcal{A}(P,\mathcal{J},I)$ hold. Assume also that $p$
  has only nondegenerate critical values in $I$, and that any two
  critical points with the same singularity type cannot have the same
  image by $p$. Then from the knowledge of
  $\Sigma(P,\mathcal{J},I)+\mathcal{O}(\h^2)$ one can infer
  \begin{enumerate}
  \item whether $E_0$ is a critical value of $p$ or not;
  \item in case $E_0$ is a critical value, the type of the
    singularity.
  \end{enumerate}
\end{theo}
This theorem is a corollary of the following proposition, where we
consider the \emph{density of states} in small regions around
$E_0$. We could, equivalently, invoke Weyl's asymptotics.
\begin{prop}
  \label{prop:density}
  Let $\gamma\in(0,1)$ and, for $E\in I$,
  \[
  \rho_\h(E) = \h^{1-\gamma}\#(\Sigma_{\h}(P,B(E,\h^\gamma))).
  \]
  Then for any $E\in I$, the limit $\rho(E)=\lim_{\h\fleche 0}
  \rho_\h(E)$ exists (in $[0,+\infty]$), and
  \begin{enumerate}
  \item if $E$ is a regular value of $p$, then $\rho$ is smooth at
    $E$;
  \item if $E$ is an elliptic critical value of $p$, then $\rho$ is
    discontinuous;
  \item if $E$ is a hyperbolic critical value of $p$, then
    $\rho(E)=+\infty$.
  \end{enumerate}
\end{prop}
\begin{demo}
  In case $E$ is a regular value, the result follows directly from
  Weyl's asymptotics, which in turn can be derived from a
  semiclassical trace formula as in~\cite{combescure-ralson-robert},
  or from the semiclassical Bohr-Sommerfeld rules as
  in~\cite{san-focus}. Let us recall the Bohr-Sommerfeld
  approach. There exists an $\epsilon>0$ such that the eigenvalues of
  $P$ inside $[E-\epsilon, E+\epsilon]$ modulo $\ohb$ are the union
  (with multiplicities) of a finite number of spectra $\sigma_k$,
  $k=1,\dots, N$, where $N$ is the number of connected components of
  $p^{-1}(E)$, and each $\sigma_k$ is determined by quasimodes
  microlocalised on the corresponding component. Precisely, the
  elements of $\sigma_k$ are given by the solutions $\lambda$ to the
  equation
  \begin{equation}
    g^{(k)}(\lambda;\h) \in 2\pi\h\ZM,
    \label{equ:bs}
  \end{equation}
  where the function $g^{(k)}$ admits an asymptotic expansion of the
  form
  \begin{equation}
    g^{(k)}(\lambda;\h) \simeq g_0^{(k)}(\lambda) + \h
    g_1^{(k)}(\lambda) + \h^2 g_2^{(k)}(\lambda) + \cdots
    \label{equ:das}
  \end{equation}
  with smooth coefficients $g_j$. Moreover, if we denote by
  $\mathcal{C}_k(\lambda)$ the $k$-ieth connected component of
  $p^{-1}(\lambda)$, in such a way that the family
  $(\mathcal{C}_k(\lambda))$ is smooth in the variable $\lambda$, then
  $g_0^{(k)}$ is the \emph{action integral}~:
  \begin{equation}
    g_0^{(k)}(\lambda) = \int_{\mathcal{C}_k(\lambda)}\xi dx.
    \label{equ:action}
  \end{equation}
  From~(\ref{equ:bs}) is follows that, for $\h$ small enough
  \[
  \#\left(\sigma_k\cap B(E,\epsilon)\right) =
  (2\pi\h)^{-1}\abs{g^{(k)}(E+\epsilon;\h)- g^{(k)}(E-\epsilon;\h)} +
  \delta,
  \]
  where the $\delta\in[-1,1]$ is here to take care of the appropriate
  integer part of the right-hand-side. Hence
  \[
  \#\left(\sigma_k\cap B(E,\epsilon)\right) =
  (2\pi\h)^{-1}\abs{2\epsilon\deriv{g_0^{(k)}(E)}{E} +
    \mathcal{O}(\epsilon^2) + \mathcal{O}(\h)} + \delta.
  \]
  With $\epsilon=\h^\gamma$, this gives
  \[
  \#\left(\sigma_k\cap B(E,\h^\gamma)\right) =
  \frac{\h^{\gamma-1}}{\pi}\abs{\deriv{g_0^{(k)}(E)}{E}} +
  \mathcal{O}(\h^{2\gamma-1}) + \mathcal{O}(1).
  \]
  Summing up all contributions for $k=1,\dots, N$, we get the first
  claim of the theorem, with
  \[
  \rho(E)= \frac1{\pi}\abs{\deriv{g_0^{(k)}(E)}{E}}.
  \]

  The second claim can be proved in a similar way, using
  Bohr-Sommerfeld rules for elliptic
  singularities~\cite{san-panoramas}. For our purposes, a Birkhoff
  normal form as in~\cite{san-charles} would even be enough, since we
  deal with energy intervals of size $\mathcal{O}(\h^\gamma)$. Here
  again there exists an $\epsilon>0$ such that the eigenvalues of $P$
  inside $[E-\epsilon, E+\epsilon]$ modulo $\ohb$ are the union (with
  multiplicities) of a finite number of spectra $\sigma_k$
  corresponding to the various connected components of
  $p^{-1}(E)$. The difference is that not all components need have
  critical points. In fact by assumption only one component may have
  an elliptic critical point. Let us call $\sigma_k$ the corresponding
  spectrum, and $\mathcal{C}_k(\lambda)$ the corresponding family of
  connected components. Since an elliptic critical point is a local
  extremum for $p$, the sets $\mathcal{C}_k(\lambda)$ are empty for
  all $\lambda$ in one of the halves of the interval $[E-\epsilon,
  E+\epsilon]$. Without loss of generality, one can assume that
  $\mathcal{C}_k(\lambda)=\emptyset$, $\forall \lambda\in
  [E-\epsilon,E[$. Then $\mathcal{C}_k(E)$ is just a point, while
  $\mathcal{C}_k(\lambda)$ is a circle for all
  $\lambda\in\,]E,E+\epsilon]$.

  The Bohr-Sommerfeld rules for elliptic singularities say that the
  elements of $\sigma_k$ are the solutions $\lambda$ to an equation of
  the form
  \begin{equation}
    e^{(k)}(\lambda;\h) \in 2\pi\h\NM,
    \label{equ:bs-ell}
  \end{equation}
  where the function $e^{(k)}$ admits an asymptotic expansion exactly
  as $g^{(k)}$ above~(\ref{equ:das}). What's more, it is equally true
  that the principal term is an action integral~:
  \[
  e^{(k)}(E)=0, \qquad e_0^{(k)}(\lambda) =
  \int_{\mathcal{C}_k(\lambda)}\xi dx, \quad \forall
  \lambda\in[E,E+\epsilon].
  \]
\end{demo}
Calculating along the same lines as above, we find, for the quantity
\[
\rho^{(k)}_\h(\lambda) := \h^{1-\gamma}\#(\sigma_k\cap
B(\lambda,\h^\gamma))),
\]
the following limits~:
\begin{enumerate}
\item when $\lambda\in[E-\epsilon,E[$, $\lim_{\h\fleche 0}
  \rho^{(k)}_\h(\lambda)=0$;
\item when $\lambda\in\,]E,E+\epsilon,E]$, $\lim_{\h\fleche 0}
  \rho^{(k)}_\h(\lambda)=
  \frac1{\pi}\abs{\deriv{e_0^{(k)}(\lambda)}{\lambda}}$;
\item $\lim_{\h\fleche 0} \rho^{(k)}_\h(E)=
  \frac1{2\pi}\abs{\deriv{e_0^{(k)}(E)}{E}}$;
\end{enumerate}

Finally, let $E$ be a hyperbolic critical value for $p$. Weyl
asymptotics for such a situation have been worked out
in~\cite{brummelhuis}, and the singular Bohr-Sommerfeld rules have
been established in~\cite{colin-p3}. Using the latter result it can be
proven as in~\cite{lablee-these} that the number of semiclassical
eigenvalues generated by a hyperbolic fixed point, in a neighbourhood
of size $\epsilon=\h^{\gamma}$ of the critical value, is of order
$\epsilon\abs{\ln\h}/\h$. Therefore, since there may be only one
hyperbolic point in $p^{-1}(0)$, it follows from this estimate and the
results we just proved above for the regular and the elliptic case
that
\[
\rho_\h(E) \geq C\abs{\ln\h},
\]
for some constant $C>0$. This gives $\rho(E)=+\infty$.
\begin{rema}
  It is probable that the nondegeneracy condition can be avoided. It
  is known quite generally that Weyl asymptotics hold for critical
  energies~\cite{zielinski-sharp}. Thus, in all case, we recover the
  action integral as the integrated density of states. It would remain
  to show that the behaviour of the action integral completely
  determines the singularities of $p$. This is easy in the Schrödinger
  case $p=\xi^2+V(x)$.
\end{rema}

\section{Topology}
\label{sec:topo}

As we already mentioned above, once the singular fibres of $p$ have
been excluded, the topology is easy to understand. The map $p$ become
a locally trivial fibration whose fibres are disjoint unions of
circles.

Thus, if $E_0$ is a regular value of $p$, the semiglobal problem
around $E_0$ just amounts to counting the number of connected
components of $p^{-1}(E_0)$.

The topology of singular fibres strongly depends on the type of
singularity. Under the nondegeneracy assumption, the topology of the
singular foliation in a neighbourhood of a singular fibre is
essentially determined by the type of the singularity, and thus by
Theorem~\ref{theo:sing}.

\subsection{Connected components}
\label{sec:connected}

Let $I$ be a compact interval of regular values of $p$. As above, we
denote by $\mathcal{C}_k(\lambda)$, for $k=1,\dots,N$ and $\lambda\in
I$ the smooth families of connected components of
$p^{-1}(\lambda)$. Each $\mathcal{C}_k(\lambda)$ is globally invariant
by the hamiltonian flow generated by $p$. Thus, this flow is periodic
on $\mathcal{C}_k(\lambda)$. Let $\abs{\tau_k(\lambda)}\neq 0$ be its
primitive period (the sign is determined by the formula below). It
follows from the action-angle theorem that $\tau_k$ is a smooth
function of $\lambda$. In fact it is well know that the period is the
derivative of the action, and we have already met this quantity in the
proof of Proposition~\ref{prop:density}. Using the action
integral~(\ref{equ:action}), we get
\[
\tau_k(\lambda)={\deriv{g_0^{(k)}(\lambda)}{\lambda}}.
\]
Notice again that $\tau_k$ never vanishes on $I$.
\begin{defi}
  We say that a point $(\lambda,t)\in(I\times \RM^*)$ is
  \textbf{resonant} whenever there exist $(k,j)$ and $(k',j')$ in
  $\{1,\dots,N\}\times\ZM^*$, with $k\neq k'$, such that
  \[
  j\tau_k(\lambda) = j'\tau_{k'}(\lambda) = -t.
  \]
\end{defi}

\begin{theo}
  \label{theo:number}
  Let $I$ be an interval of regular values of $p$, and let Assumption
  $\mathcal{A}(P,\mathcal{J},I)$ hold.  Assume also that the set of
  resonant points in $I\times\RM$ is discrete.  Then the number $N$ of
  connected components of $p^{-1}(\lambda)$, $\lambda\in I$, is
  determined by spectrum $\Sigma(P,\mathcal{J},I)+\mathcal{O}(\h^2)$.
\end{theo}
Before proving the theorem, let us just remark that the leading term
of Weyl's asymptotics is not sharp enough for this. Indeed, it only
gives the density $\rho$ (Proposition~\ref{prop:density})~:
\begin{equation}
  \rho(\lambda) = \frac1{\pi}\sum_{k=1}^N \abs{\tau_k(\lambda)}.
  \label{equ:density}  
\end{equation}
From this one cannot distinguish, for example, one component with
period $\tau$ from two components with periods
$\abs{\tau_1}+\abs{\tau_2}=\abs{\tau}$.

Remark also that the condition on resonant points is not adapted to
systems with symmetries. For instance, a Schrödinger operator with a
symmetric double well has two components with equal periods.

\vspace{2ex}

\begin{demo}[of Theorem~\ref{theo:number}]
  We introduce the period lattice $\mathcal{L}_k(I)$~:
  \[
  \begin{split}
    \mathcal{L}_k(I) :&=  \{(\lambda,t)\in I\times \RM;\quad
    \exp(t\ham{p}) \text{
      is periodic on } \mathcal{C}_k(\lambda)\}\\
    & = \{(\lambda,j\tau_k(\lambda)); \quad \lambda\in I, j\in\ZM\}
    \label{equ:lattice}
  \end{split}
  \]
  and $\mathcal{L}(I)=\bigcup_{k=1}^N\mathcal{L}_k(I)$.  The set
  $\mathcal{L}(I)$ is a union of smooth graphs that may intersect. The
  intersection points for $t\neq 0$ are precisely the resonant points.

  In order to prove the theorem, we split the argument into two
  steps. The first one is to prove that
  $\Sigma(P,\mathcal{J},I)+\mathcal{O}(\h^2)$ determines
  $\mathcal{L}(I)$. The second step consists in showing why the
  knowledge of $\mathcal{L}(I)$ --- and the assumption on the set of
  resonant points --- allows us to count the number $N$ of connected
  components.

\paragraph{Step 1.}
Coming back to the Bohr-Sommerfeld rules discussed in the proof of
Proposition~\ref{prop:density}, we recall that the spectrum of $P$
modulo $\ohb$ is the superposition of the spectra $\sigma_k$ generated
by $\mathcal{C}_k$, for $k=2,\dots,N$. For each $k$, $\sigma_k$ has a
periodic structure that makes it close to an arithmetic
progression. Thus, a simple and naive idea to distinguish between the
different periodic structures is to perform a \emph{frequency
  analysis}, via a Fourier transform. Because we have at our disposal
only a truncated sequence of eigenvalues (those that belong to $I$),
we need to introduce a cut-off. Let $I'\Subset I$ and let
$\chi\in\Cinf(\RM)$ have compact support in the interior of $I$ and be
equal to $1$ on $I'$. We introduce the spectral measure
\[
D_0(\lambda;\h) = \sum_{E\in\Sigma_\h(P,I)} \chi(E) \delta_E(\lambda),
\]
where $\delta_E$ is the Dirac distribution at $E$. The quantity we
want to investigate is its Fourier transform. Since the mean spacing
between consecutive eigenvalues is of order $\h$, we use a
corresponding scale for the time variable $t$, and thus introduce
\[
Z(t;\h)=\sum_{E\in\Sigma_\h(P,I)} \chi(E)e^{-itE/\h}.
\]
The function $Z$ is called the \emph{partition function}. In fact, the
idea we've just described is very well known in the semiclassical
context, and is part of the general formalism of \emph{trace
  formul\ae}. We can consider the Schrödinger group
$U(t;\h)=\exp(-itP/\h)$, and then
$Z(t;\h)=\textup{Trace}(\chi(P)U(t;\h))$. It is well known that
$\chi(P)U(t;\h)$ is a Fourier Integral Operator, whose canonical
transformation is the classical flow of $p$. Moreover, its trace is a
lagrangian (or WKB) distribution associated with the lagrangian
manifold of periods
\[
\Lambda_p=\{(E,\tau)\in\RM^2; \quad \exists z\in p^{-1}(E),
\exp(\tau\ham{p})(z)=z\} = \mathcal{L}(I').
\]
Such a result would almost finish the proof of Step 1. In fact, this
statement exists in many versions, depending on various possible
situations and hypothesis.  For this reason we are not using it here
as is, but instead resort once again to the Bohr-Sommerfeld rules,
which is arguably the easiest way to go.

We can split the partition function as
\[
Z(t;\h) = \sum_{k=1}^N\sum_{E\in\sigma_k} \chi(E)e^{-itE/\h}.
\]
Then from~(\ref{equ:bs}) one can introduce $c\mapsto f^{(k)}(c;\h)$ as
the inverse of $\lambda\mapsto g^{(k)}(\lambda;\h)$, which exists for
$\h$ small enough, and write
\begin{equation}
  Z(t;\h) = \sum_{j\in\ZM} \phy_t(2\pi\h j;\h)
  \label{equ:partition}
\end{equation}
(which, as before, is a finite sum) with
\begin{equation}
  \phy_t(c;\h):=\sum_{k=1}^N \chi(f^{(k)}(c;\h))
  e^{-itf^{(k)}(c;\h)/\h}.
  \label{equ:phy}
\end{equation}
Note that $\phy_t(\cdot;\h)\in\Cinf_0(\RM)$. By the Poisson summation
formula,
\begin{equation}
  Z(t;\h) = \frac{1}{2\pi\h}\sum_{j\in\ZM} \hat{\phy}_t(j/\h)
  \label{equ:Z}
\end{equation}
(which, contrary to~(\ref{equ:partition}), is a truly infinite sum)
with
\begin{equation}
  \hat{\phy}_t(j/\h) = \int_{\RM} e^{-icj/\h}\phy_t(c)dc =
  \sum_{k=1}^N Z_k(t;j,\h)\label{equ:phyhat}
\end{equation}
and
\[
Z_k(t;j,\h) = \int e^{-i\h^{-1}(cj +
  tf^{(k)}(c;\h))}\chi(f^{(k)}(c;\h)) dc.
\]
The integral $Z_k$ is a compactly supported oscillatory integral,
whose phase is stationary when $j+t\deriv{f_0^{(k)}}{c} = 0$ or
equivalently
\begin{equation}
  t= -j\tau_k(\lambda), \qquad \lambda=f^{(k)}(c;\h).\label{equ:phase}
\end{equation}
Moreover, the Hessian of the phase, $Q:=t\frac{\partial^2
  f_0^{(k)}}{\partial c^2}=t(\deriv{\tau_k}{\lambda})^{-1}$ never
vanishes for $t\neq 0$. Hence, by the stationary phase expansion,
$Z_k$ is a lagrangian distribution whose principal symbol $\iota_k$
can be written as a smooth function of $\lambda$~:
\[
\iota_k(\lambda;j,\h)=
\frac{e^{i\frac{\pi}4\text{sign}(Q)}}{\abs{Q}^{1/2}}
e^{-ij\h^{-1}(g_0^{(k)}(\lambda)-\lambda\tau_k(\lambda))}\chi(\lambda),
\]
with $Q=-j\tau_k(\lambda)(\deriv{\tau_k}{\lambda})^{-1}$.  Since its
amplitude vanishes precisely with $\chi$, we can deduce that the
semiclassical wave-front of $Z_k$ is (for fixed $j\in\ZM^*$)
\[
WF_{\h}(Z_k) = \{(\lambda,t)\in\RM^2; \quad t=-j\tau_k(\lambda),
\chi(\lambda)\neq 0\}.
\]
We still need to sum up all $Z_k(t;j,\h)$ for $j\in\ZM^*$.
For this we consider the localisation of $Z$. Without loss of
generality, one can restrict to positive times.  Let $t_0> 0$,
$\epsilon>0$, and let $\rho\in\Cinf_0(B(t_0,\epsilon))$. There is no
solution to~(\ref{equ:phase}) in the support of $\rho$ for $\abs{j}$
outside the interval
\[
I_k(\epsilon):=\left(\frac{t_0-\epsilon}{\sup_J\abs{\tau_k}},
  \frac{t_0+\epsilon}{\inf_J\abs{\tau_k}}\right).
\]
Making explicit the non-stationary phase argument, we can write, for
any $\ell\in\NM$,
\[
Z_k(t;j,\h) =\left(\frac{\h}{ji}\right)^\ell \int e^{-i\h^{-1}(cj +
  tf_0^{(k)}(c;\h))} L^\ell(a(c;\h)) dc,
\]
where $L$ is the linear differential operator defined by
\[
(Lu)(c) =
\frac{d}{dc}\left(\frac{u(c)}{1+\frac{t}{j}\deriv{f_0^{(k)}}{c}}\right)
\]
and $a(\cdot;\h)\in\Cinf_0(I)$ admits an asymptotic expansion in
non-negative powers of $\h$, in the $\Cinf$ topology. Let
$b(c)=(1+\frac{t}{j}\deriv{f_0^{(k)}}{c})^{-1}$. Then $b$ is uniformly
bounded on $I$ for $\abs{j}>(t_0+\epsilon)/\inf_J\abs{\tau_k}$, and
for any $\ell\in\NM^*$, there exists a positive constant $C_{\ell}$,
independent of $j$ and $\h$, such that
$\abs{\frac{d^{\ell}b}{dc^{\ell}}}\leq C_{\ell}/j$. Therefore, there
exist constants $\tilde{C}_{\ell}>0$ such that
\[
\abs{L^{\ell}(a)} \leq \tilde{C}_{\ell},
\]
and we get, again when $\abs{j}>(t_0+\epsilon)/\inf_J\abs{\tau_k}$,
\[
\abs{\rho(t)Z_k(t;j,\h)} \leq \tilde{C}_{\ell}
\left(\frac{\h}{j}\right)^\ell.
\]
Thus, for $\ell\geq 2$,
\[
\sum_{\abs{j}>\frac{t_0+\epsilon}{\inf_J\abs{\tau_k}}}
\abs{\rho(t)Z_k(t;j,\h)} \leq \tilde{C}_\ell \h^{\ell}.
\]
This shows that only a finite (independent of $\h$) number of terms
contribute to $\rho(t)Z(t;\h)$ modulo $\ohb$. Thus the
(non)-stationary phase approximations are jointly valid. Therefore
$Z(t;\h)$ microlocally vanishes at any point that does not belong to
$\mathcal{L}(I)$; this writes
\[
WF_{\h}(Z(\cdot;\h)) \subset \mathcal{L}(I).
\]
More precisely,
\[
WF_{\h}(\rho Z(\cdot;\h)) \subset \left\{(\lambda,j\tau_k(\lambda));
  \quad \lambda\in I, \abs{j} \in I_k(\epsilon), k=1,\dots,N\right\}.
\]
Moreover, at a non-resonant point $(\lambda,j\tau_k(\lambda))$, no
other period $j'\tau_{k'}$ can contribute, and $Z(\cdot;\h)$ is a
lagrangian distribution with principal symbol equal to
$\iota_k(\lambda;j,\h)$. Since the set of resonant points is discrete,
and $WF_{\h}(Z)$ is closed in $T^*I$, we must have
$WF_{\h}(Z(\cdot;\h)) = \mathcal{L}(I)$, which finishes the proof of
the first step.

\paragraph{Step 2.} We are now left with a simple geometric inverse
problem~: given the set of periods $\mathcal{L}(I)$, how can one
recover the number $N$ of connected components ?

Our strategy is to recover the fundamental periods
$\abs{\tau_1},\dots,\abs{\tau_N}$. First of all, by Weyl's
asymptotics~(\ref{equ:density}), one obtains the \emph{a priori} bound
$\abs{\tau_k(\lambda)}\leq \pi\rho(\lambda)$. Let $R:=\max_J
\pi\rho$. Then by assumption, the set of resonant points inside
${I}\times ]0,R]$ is finite; therefore, one can always find a smaller,
non-empty interval $\tilde{I}\subset I$ such that there is \emph{no}
resonant point at all in $\tilde{I}\times ]0,R]$.

We extract the periods $\tau_k$ from
$\mathcal{L}_1:=\mathcal{L}(\tilde{I})\cap (\tilde{I}\times ]0,R])$
inductively, as follows.

\begin{enumerate}
\item Consider a point $(\lambda_1,\tau_1)\in\mathcal{L}_1$ with
  ``minimal height'' $\tau_1$~: $\forall
  (\lambda,\tau)\in\mathcal{L}_1, \tau_1\leq \tau$.
\item By the non-resonance assumption, the connected component of
  $(\lambda_1,\tau_1)$ in $\mathcal{L}_1$ is the graph of a smooth
  function of the interval $\tilde{I}$. We denote this function by
  $\lambda\mapsto\tau_1(\lambda)$.
\item Consider the set
  \[
  \mathcal{L}_2:=\mathcal{L}_1\setminus\{(\lambda,j\tau_1(\lambda));
  \quad \lambda\in\tilde{I}, j\in\ZM^*\}.
  \]
  Again by the non-resonance assumption, $\mathcal{L}_1$ remains a
  union of non-intersecting smooth graphs.
\item If $\mathcal{L}_1$ is empty, then $N=1$. Otherwise, start again
  by replacing $\mathcal{L}_0$ by $\mathcal{L}_1$, and so on.  If
  $\mathcal{L}_k$ is empty, then $N=k-1$.
\end{enumerate}
\end{demo}

\begin{rema}
  If we disregard symmetry issues, our assumption on the resonant set
  is quite weak. For instance, one can allow the crossing of two
  periods to be flat (all derivatives are equal at a point $\lambda$),
  simply because we put ourselves in a region with no crossing at
  all. However, it is easy to prove Step 2 with even weaker
  assumptions. For instance, it may work even if there are some open
  intervals of values of $\lambda$ which admits resonant pairs.  It
  would be interesting to know whether Step~1 could hold in this case
  as well. It would then involve sub-principal terms in the
  Bohr-Sommerfeld expansion.
\end{rema}

\subsection{Singular fibres}

As we already mentioned, the following result comes for free.
\begin{theo}
  \label{theo:singular_fibres}
  Let Assumption $\mathcal{A}(P,\mathcal{J},I)$ hold, and
  let $E_0\in I$ be a nondegenerate critical value of $p$.  Assume
  also that $p^{-1}(E_0)$ contains only one critical point.  Then from
  the knowledge of $\Sigma(P,\mathcal{J},I)+\mathcal{O}(\h^2)$ one can
  determine the topology of the singular foliation induced by $p$, in
  a saturated neighbourhood of $p^{-1}(E_0)$.
\end{theo}
\begin{demo}
  Under these assumptions, the topology of the singular foliation
  induced by $p$ in a saturated neighbourhood of $p^{-1}(E_0)$ is known
  to be completely characterised by the type of the
  singularity~\cite{dufour-mol-toul,zung-I}, which is determined by
  Theorem~\ref{theo:sing}. For the convenience of the reader, we
  briefly recall the two possible cases.
  \begin{enumerate}
  \item \emph{The elliptic case. ---} The singular fibre $p^{-1}(E_0)$
    is just a point and the foliation is homeomorphic to the one given
    by the Hamiltonian $H(x,\xi)=x^2+\xi^2$.
  \item \emph{The hyperbolic case. ---} The singular fibre is a circle
    with a transversal self-intersection (the figure eight). It separates
    a saturated neighbourhood into three connected parts~: two on one
    side, and one on the other side. It is homeomorphic to the
    foliation given by the Hamiltonian $H(x,\xi)=\xi^2+x^4-x^2$, in a
    neighbourhood of $H^{-1}(0)$.
  \end{enumerate}
\end{demo}

\subsection{Global topology}

We say that a hamiltonian system $p$ on the symplectic 2-manifold $M$
is topologically equivalent to the hamiltonian system $\tilde{p}$ on
$\tilde{M}$ if there is a homeomorphism $\phy:M\fleche \tilde{M}$ such
that
\[
p = \tilde{p}\circ\phy.
\]

Notice that this implies that $\phy$ respects the foliation, fibre by
fibre. In particular, $p$ and $\tilde{p}$ have the same set of regular
values and the same set of critical values. If $I$ is an open
interval, then two hamiltonian systems $p$ and $\tilde{p}$ are called
topologically equivalent over $I$ when they are topologically
equivalent when restricted to the symplectic manifolds $p^{-1}(I)$,
$(\tilde{p})^{-1}(I)$.

We call the \emph{topological type} of a hamiltonian system the
equivalence class of topologically equivalent systems.

\begin{theo}
  \label{theo:topo}
  Let Assumption $\mathcal{A}(P,\mathcal{J},I)$ hold, and assume that
  $p$ has only nondegenerate critical values in some neighbourhood of
  $I$, such that any two critical points with the same singularity
  type cannot have the same image by $p$. Let $c_1<\cdots<c_n$ be the
  critical values of $p$ in $I$. Suppose that in each interval
  $(c_i,c_{i+1})$, $i=1,\dots,n-1$, there exists a non-empty
  subinterval $I_i$ such that the set of resonant points in
  $I_i\times\RM$ is discrete in $\RM^2$. Then the knowledge of
  $\Sigma(P,\mathcal{J},I)+\mathcal{O}(\h^2)$ determines the
  topological type of the hamiltonian system $p$ over $I$.
\end{theo}
\begin{demo}
  Upon a possible enlargement of $I$, one may assume that
  $I=(E_0,E_1)$ for regular values $E_0,E_1$. Using symplectic
  cutting~\cite{lerman-cut} or surgery~\cite{zung-I}, one may replace
  the phase space $\RM^4$ by a compact symplectic manifold where
  $p^{-1}(I)$ is embedded. Then we apply the result
  of~\cite{dufour-mol-toul} that says that the topological type of $p$
  on $M$ is determined by its Reeb graph~: the set of leaves of the
  foliation, as a topological 1-complex. This graph is characterised
  by the relative positions of critical values, and the number of
  fibres between two consecutive critical values. The former is
  determined by the spectrum in $I$ thanks to Theorem~\ref{theo:sing},
  while the latter is determined for each $i=1,\dots,n-1$ by the
  spectrum in $I_i$, thanks to Theorem~\ref{theo:number}. This give
  the topological type of $p$, up to some homeomorphism of the Reeb
  graph itself. But since we know the precise values of $p$ at
  singularities, we can in fact assume that this homeomorphism is the
  identity.
\end{demo}

\section{Symplectic geometry}
\label{sec:sympl}

The tools we've used so far give us the \emph{periods} of the
classical hamiltonian system, which is of course much more than a mere
topological information. We show here that it is indeed sufficient to
recover the full dynamics of the systems.

We say that a hamiltonian system $p$ on the symplectic 2-manifold $M$
is symplectically equivalent to the hamiltonian system $\tilde{p}$ on
$\tilde{M}$ if there is a smooth symplectomorphism $\phy:M\fleche
\tilde{M}$ such that
\[
p = \tilde{p}\circ\phy.
\]

Thus, the dynamics of $p$ on the levelset $\{p=E\}$ is transported via
$\phy$ to the dynamics of $\tilde{p}$ on the levelset
$\{\tilde{p}=E\}$.

We call the \emph{symplectomorphism type} of a hamiltonian system the
equivalence class of symplectically equivalent systems. As before, one
may restrict this equivalence to an interval $I$ of values of $p$ and
$\tilde{p}$.

\begin{defi}
  Let $(\lambda,t)\in I\times \RM$ be a resonant point for $p$. Thus
  \[
  j\tau_k(\lambda) = j'\tau_{k'}(\lambda) = -t
  \]
  for some $j,j',k\neq k'$. We say that this resonance is \emph{weakly
    transversal} if there exists an integer $n\in\NM^*$ such that the
  $n$-th derivatives of the periods are not equal~:
  \[
  j\tau_k^{(n)}(\lambda) \neq j'\tau_{k'}^{(n)}(\lambda).
  \]
\end{defi}
\begin{theo}
  \label{theo:sympl}
  Let Assumption $\mathcal{A}(P,\mathcal{J},I)$ hold, and suppose that
  $p$ has only nondegenerate critical values in some neighbourhood of
  $I$, such that any two critical points with the same singularity
  type cannot have the same image by $p$. Let $c_1<\cdots<c_n$ be the
  critical values of $p$ in $I$. Suppose that for each interval
  $J_i:=(c_i,c_{i+1})$, $i=1,\dots,n-1$, the set of resonant points in
  $J_i\times\RM$ is discrete. Finally assume that all such resonant
  points are weakly transversal.

  Then the knowledge of $\Sigma(P,\mathcal{J},I)+\mathcal{O}(\h^2)$
  determines the symplectic type of the hamiltonian system $p$ over
  $I$.
\end{theo}
\begin{demo}
  We use the symplectic classification
  of~\cite{dufour-mol-toul,toulet-these} using \emph{weighted} Reeb
  graphs. Under our assumptions, the Reeb graph has vertices of degree
  1 and 3. A vertex of degree 1, a \emph{bout}, corresponds to an
  elliptic critical value, while a vertex of degree 3, called a
  \emph{bifurcation point}, corresponds to a hyperbolic critical
  value.  At a bifurcation point we can distinguish one particular
  edge, called the \emph{trunk}, corresponding to the side of the
  figure 8 with only one connected component. The two other edges are
  called the \emph{branches}.  A weighted Reeb graph is a Reeb graph
  each of whose edges is associated with a positive real number, its
  \emph{length}, and such that each of the two branches of each
  bifurcation point is associated with a formal Taylor series
  (\emph{ie} a sequence of real numbers). The hypothesis of the
  theorem allow for determining the topological Reeb graph via
  Theorem~\ref{theo:topo}.  Thus, the next step of the proof is to
  show how the numbers that constitute the weighted Reeb graph can be
  recovered from the spectrum. The final step is to obtain the
  symplectic equivalence in the sense that we have just defined above.

\paragraph{The lengths. ---}
Let $\mathcal{C}_k(J_i)$, for $k=1,\dots,N$, be the connected
components of $p^{-1}(J_i)$. Let $K_{k,i}\in\Cinf(\mathcal{C}_k(J_i))$
be an action variable for the regular lagrangian fibration $p_{\restr
  \mathcal{C}_k(J_i)}$; it is unique up to a sign and an additive
constant. By definition the \emph{length} of the edge corresponding to
the set of leaves in $\mathcal{C}_k(J_i)$ is
\begin{equation}
  \ell_{k,i}:=\abs{\lim_{c\fleche c_i}K_{k,i}(c) - \lim_{c\fleche
      c_{i+1}}K_{k,i}(c)}.
  \label{equ:length}
\end{equation}
In a learned terminology, this is the Duistermaat-Heckman measure of
$J_i$ for the $S^1$-action defined by $K_{k,i}$, or, equivalently, it
is the affine length of $J_i$ endowed with its natural integral affine
structure given by $p_{\restr \mathcal{C}_k(J_i)}$.

It follows from the local models for elliptic and hyperbolic
singularities that this length is always finite. This is obvious at
elliptic singularities, where the action has the form $x^2+\xi^2$. At
a hyperbolic singularity $m$, one can introduce a foliation function
$q$ such that, in some local symplectic coordinates around $m$,
$q=x\xi$, and $q>0$ on the branches while $q<0$ on the trunk. Then the
Duistermaat-Heckman measure has the form
\begin{equation}
  \begin{cases}
    d\mu_j(q) & = \left(\ln q + g_j(q)\right)dq \quad \text{ on each
      branch (j=1,2)}\\
    d\mu(q) & = \left(2\ln \abs{q} + g(q)\right)dq \quad \text{ on the
      trunk},
  \end{cases}
  \label{equ:hyperbolic}
\end{equation}
with some smooth functions $g,g_1,g_2$ satisfying
\[
\forall p, \qquad g^{(p)}(0) = g_1^{(p)}(0)+ g_2^{(p)}(0).
\]
Under this form, the Taylor series of the functions $g$, $g_1$, $g_2$
at the origin are uniquely
defined~\cite{toulet-these,san-semi-global}.

Using the proof of Theorem~\ref{theo:number}, from the spectrum in $I$
we can recover the periods $\tau_k(\lambda)$, $k=1,\dots,N$, for
$\lambda$ in any interval in $J_i$ where the graphs of the periods
$\tau_k$ don't cross. At a crossing the difficulty is to put the
labels $k$ correctly, so that the connected components
$\mathcal{C}_k(\lambda)$ remain in the same $\mathcal{C}_k(J_i)$ when
$\lambda$ varies.  This can be overcome precisely thanks to the weak
resonant assumption at each crossing, because each $\tau_k$ is $\Cinf$
in $J_i$. This was the main issue. Now, fixing a point $\lambda_i\in
J_i$, the action variable $K_{k,i}$ can be computed by the formula
\[
K_{k,i}(\lambda):=\int_{\lambda_i}^\lambda \tau_k(\lambda)d\lambda,
\quad \lambda\in J_i.
\]
This gives the length of $\mathcal{C}_k(J_i)$ via
equation~(\ref{equ:length}).

\paragraph{The Taylor series at the bifurcation points. ---}
By definition, the sequences of numbers associated with a bifurcation
point in the Reeb graph are the Taylor series of the functions $g_1$,
$g_2$ (defined in equation~(\ref{equ:hyperbolic})) at the origin.

Let us show how to recover the Taylor series of $g$ from the
spectrum. The procedure is completely analogous for $g_1$ and $g_2$.

Thus, we consider a hyperbolic critical value $c_{i+1}$.  We want to
express the Duistermaat-Heckman measure on the trunk in terms of the
principal symbol $p$. By a theorem of Colin de Verdière and
Vey~\cite{colin-vey}, there exist local symplectic coordinates
$(x,\xi)$ at the hyperbolic point, and a smooth, locally invertible
function $f:(\RM,c_{i+1})\fleche (\RM,0)$ such that
\[
f(p)=x\xi=q.
\]
For notational purposes, one may assume that $f'(c_{i+1})>0$, which
amounts to say that the trunk is sent by $p$ to $\lambda<c_{i+1}$.
Then from~(\ref{equ:hyperbolic}), for $\lambda$ close to $c_{i+1}$,
$\lambda<c_{i+1}$,
\[
d\mu(\lambda) = \left(2\ln\abs{f(\lambda)} + g\circ f(\lambda)\right)
f'(\lambda)d\lambda.
\]
On the other hand if the connected component corresponding to the
trunk is $\mathcal{C}_k(J_i)$, one has by definition of the
Duistermaat-Heckman measure
$d\mu(\lambda)=\tau_k(\lambda)d\lambda$. Therefore
\[
\tau_k(\lambda) = f'(\lambda)\left(2\ln\abs{f(\lambda)} + g\circ
  f(\lambda)\right) = 2f'(\lambda)\ln\abs{\lambda-c_{i+1}} +
h(\lambda),
\]
for some smooth function $h$ at $\lambda=c_{i+1}$. There, using
Taylor's formula, we have written $f(\lambda)=\alpha(\lambda-c_{i+1})
+ (\lambda-c_{i+1})^2\hat{f}(\lambda)$, with $\alpha>0$ and $\hat{f}$
smooth at $c_{i+1}$, and hence
\begin{equation}
  h(\lambda) =
  2f'(\lambda)\ln\abs{\alpha+(\lambda-c_{i+1})\hat{f}(\lambda)} +
  f'(\lambda)g\circ f(\lambda).
  \label{eqh:h}
\end{equation}
This shows that $h$ is smooth for $\lambda$ close to $c_{i+1}$.

It is easy to see that any smooth function $\phi$ in a neighbourhood of
the origin such that $\phi(t)\ln t$ extends to a smooth function at
$t=0$ must be flat. Hence the knowledge of $\tau_k(\lambda)$ for
$\lambda<c_{i+1}$ completely determines the Taylor series of
$f'(\lambda)$ (and hence $f(\lambda)$) at $\lambda=c_{i+1}$.

Then one can recover the Taylor series of $h$ using
\[
h(\lambda)=\tau_k(\lambda)-2f'(\lambda)\ln\abs{\lambda-c_{i+1}}, \quad
\forall\lambda<c_{i+1}
\]
Finally, from~(\ref{eqh:h}) and the fact that $f$ is locally
invertible, one can recover the Taylor series of $g$ at the origin.

\paragraph{Symplectic equivalence. ---}
We have proven that the weighted Reeb graph is determined by the
spectrum. By Toulet's
classification~\cite{dufour-mol-toul,toulet-these}, if two such
systems $(M,p)$ and $(M,\tilde{p})$ have the same weighted Reeb graph,
there exists a symplectomorphism $\phy:M\fleche \tilde{M}$ such that
$p$ and $\tilde{p}\circ\phy$ define the same singular foliation on $M$
($\phy$ indices a homeomorphism of the leaf space, fixing the
vertices). If we assume that the operators $P$ and $\tilde{P}$ have
the same spectrum (modulo $\h^2$) and fulfil the requirements of the
theorem, then we also know that $p$ and $\tilde{p}\circ\phy$ share the
same set of critical values $c_i$. The fact that $p$ and
$\tilde{p}\circ\phy$ define the same foliation implies that for each
connected component $\mathcal{C}_k(J_i)$, there exists a smooth,
invertible function $f:J_i\fleche J_i$ such that
\begin{equation}
  p = f\circ \tilde{p} \circ \phy \qquad \text{ on } \mathcal{C}_k(J_i).
  \label{equ:equiv}
\end{equation}
Since the singular fibres at the ends of $\mathcal{C}_k(J_i)$ are
fixed by $\phy$, $f$ must be increasing, and thus extends to a
homeomorphism of $\overline{J_i}$.

As we already saw, the spectrum also determines the periods at a given
energy $E=\lambda$. Hence for $\lambda\in J_i$,
$\tau_k(\lambda)=\tilde{\tau}_k(\lambda)$. Since $\tau_k$ is
integrable at $c_{i+1}$, we can define action integrals for
$\lambda<c_{i+1}$ as~:
\[
K_{k,i}(\lambda):=\int_{c_{i+1}}^\lambda \tau_k(\lambda)d\lambda,
\qquad \tilde{K}_{k,i}(\lambda):=\int_{c_{i+1}}^\lambda
\tilde{\tau}_k(\lambda)d\lambda.
\]
We have $K_{k,i}(\lambda)=\tilde{K}_{k,i}(\lambda)$. On the other
hand, the action is a symplectic invariant of the
foliation. From~(\ref{equ:equiv}) on can compute the action on the
curve
$\phy(\mathcal{C}_k(f(\lambda)))=\tilde{\mathcal{C}}_k(\lambda)$~:
$K_{k,i}(f(\lambda))=\tilde{K}_{k,i}(\lambda)+\textup{const}$. Therefore
\[
K_{k,i}(\lambda)=K_{k,i}(f(\lambda)).
\]
Since $\tau_k$ does not vanish in $J_i$, $K_{k,i}$ is strictly
monotonous on $J_i$. Therefore
\[
f(\lambda)=\lambda, \qquad \forall \lambda\in J_i.
\]
Thus $p=\tilde{p}\circ \phy$ on each $\mathcal{C}_k$, and by
continuity
\[
p=\tilde{p}\circ\phy \quad \text{ on } M.
\]
This finishes the proof of the theorem.
\end{demo}

\bibliographystyle{plain}
\bibliography{bibli-utf8}
\end{document}